\documentclass[11pt]{article}

\usepackage{comment}
\usepackage[linktocpage=true]{hyperref}
\usepackage[margin=1.2 in, top=1 in, bottom= 1.2 in]{geometry}

\usepackage{amsfonts}
\usepackage{mathrsfs}
\usepackage{amssymb}
\usepackage{mathtools}
\usepackage{tabularx}
\usepackage{makecell}
\usepackage{arydshln}

\usepackage{enumitem}
\setlist{nolistsep}
\usepackage{amsthm}
\usepackage{relsize}
\usepackage{nicefrac}
\usepackage[capitalize]{cleveref}

\newtheoremstyle{plain}{3mm}{3mm}{\slshape}{}{\bfseries}{.}{.5em}{}
\newtheoremstyle{definition}{2mm}{2mm}{}{}{\bfseries}{.}{.5em}{}
\theoremstyle{plain}

\newtheorem{theorem}{Theorem}

\theoremstyle{definition}

\theoremstyle{plain}
\newtheorem*{namedthm}{\namedthmname}
\newcounter{namedthm}
	

\title{On certain partition bijections related to Euler's partition problem.}
\date{}
\author{Aritro Pathak}

\begin{document}

\maketitle

\begin{abstract}
     We give short elementary expositions of combinatorial proofs of some variants of Euler's partition identity that were first addressed analytically by George Andrews, and later combinatorially by others. The method using certain matrices to concisely explain these bijections, enables us to also give new generalizations of two of these results.
\end{abstract}

\section{Introduction and statement of results}
There appeared a couple of conjectures \cite{one},\cite{two}, on partition identities in the On-line Encyclopedia of Integer sequences, which were generalizations of Euler's famous partition result that the number of partitions of any positive integer $n$ into all distinct parts equals the number of partitions where each part is odd. These results were first proven using generating functions by George Andrews \cite{Andrews}.

We state the two conjectures here:

\begin{theorem}\label{thm:conj1}
 Let $n$ be a positive integer, and let $a(n)$ be the number of partitions of $n$ that contain exactly one even part, which may be repeated more than once. Then $a(n)$ is also the difference between the total number of parts in the odd partitions of $n$ and the total number of parts in the distinct partitions of $n$.
\end{theorem}

\begin{theorem}\label{thm:conj2}
Let $n$ be a positive integer, and let $a_1(n)$ be the number of partitions of $n$ such that there is exactly one part appearing three times and other parts appear exactly once. Then $a_{1}(n)$ is also the difference between the total number of parts in the distinct partitions of $n$ and the number of distinct parts in the odd partitions of $n$.
\end{theorem}

After Andrews' analytical proof of these results, combinatorial proofs of the above conjectures were given in the papers of Yang \cite{Yang}, Fu and Tang \cite{Fu}, and Ballantine and Bielak \cite{Ballantine}. Another direction in which we could study a variant of the statement of \cref{thm:conj1} is to see if we can allow more than one even part. This is what Andrews asked at the end of his paper \cite{Andrews}. Indeed, we have the following theorem, which is the $k=2$ case of Theorem 1.4 of Fu and Tang's paper \cite{Fu}.

\begin{theorem}\label{thm:conj3}

Let $n$ and $k$ be positive integers. Let $A_k(n)$ be the number of partitions of $n$ where there are exactly $k$ distinct even parts, each possibly repeated. Let $B_k(n)$ be the number of partitions that have exactly $k$ repeated parts. Then $A_k(n)=B_k(n)$.

\end{theorem}

There is a well known generalization of Euler's partition identity, which was first proven by Glaisher \cite{lehmer}. Glaisher originally used his technique to give the first known combinatorial proof of Euler's partition identity, and generalized it to the following:

\begin{theorem}\label{thm:oldone}
Given an integer $d>0$, the number of partitions $c(n)$ of an integer $n$ such that no part is divisible by $d$, is the same as the number of partitions $e(n)$ of $d$ where no part appears more than $d$ times.
\end{theorem}

Clearly the $d=2$ case is Euler's partition theorem. Recently \cite{Pathak} a generalization was given of the original proof of Glaisher, by looking at the complementary problem: showing that there is a bijective correspondence between the set of partitions of $n$ where at least one part appears $d$ times and the set of partitions where at least one part is divisible by $d$. This generated a rich family of bijections; and seems to have not appeared in the literature before. We briefly outline this proof in Section 2.


In this paper we give quick elementary expositions of combinatorial proofs of \cref{thm:conj1}, \cref{thm:conj2}, and \cref{thm:conj3} but which essentially reduce to being variants of the proofs earlier given in \cite{Yang,Fu,Ballantine}. We then state and prove two  results, one of which extends \cref{thm:conj2} and the other extends \cref{thm:conj3}. We state them below:

\begin{theorem}\label{thm:new1}
Given a fixed positive integer $p\geq 2$, consider the set of partitions $A_{k}(n)$ of any given positive integer $n$, so that there are exactly $k$ parts whose highest factor of $2$ has exponent greater than or equal to $p$. Consider the set of partitions $B_{k}(n)$ of $n$ that are such that exactly $k$ parts appear at least $2^{p}$ times. Then there is a family of bijections between these two sets of partitions.

\end{theorem}

 Note that \cref{thm:new1} is also true for the case $p=1$, which is \cref{thm:conj3}. Note that by taking $k=2^p$ for $p\geq 2$, we get this result from Theorem 1.4 of Fu and Tang's paper \cite{Fu}, however we get a larger family of bijections for the specific case of $p\geq 2$ in our theorem, and essentially the Fu Tang bijection works differently. As an example, for any $p\geq 2$ and odd integer $\alpha$, taking $k=2^p$ in Fu and Tang theorem, the elements $\alpha \cdot 2^p$ and $\alpha \cdot 2^{p+1}$ belong to separate matrices, whereas in our bijective method, we would consider both these elements to be within the same matrix $M^{(\alpha,\lambda)}$ the construction of which is outlined in the next section. In our case we easily get a rich family of bijections in \cref{thm:new1} in the spirit of the proof in \cite{Pathak}. The statement of \cref{thm:new1} also naturally holds true when we consider any arbitrary integer in place of $2$ and look at it's exponents, and here too we would get a larger family of bijections than from the earlier result.

Further we state this next new generalization that extends \cref{thm:conj3}.

\begin{theorem}\label{thm:new2}

For a positive integer $n$, let $f(n)$ be the number of partitions of $n$ such that there is exactly one part appearing five times, and all other parts appear once. Also consider the set $G(n)$ of distinct partitions of $n$ with the property that if $\alpha \cdot 2^{i} (i\geq 0 )$ appears in the partition with $\alpha$ any odd integer, then $\alpha \cdot 2^{i+1}$ does not appear in the partition. Also consider the set $H(n)$ of partitions of $n$ with only odd parts such that in the base $2$ expansion of the number of times any odd number $\beta$ appears, there are no two consecutive 1's. Then $f(n)$ is exactly the difference of the number of distinct parts appearing in $G(n)$ and the number of distinct parts in $H(n)$.
\end{theorem}

\section{Combinatorial proofs of Theorems 3 and a generalization.}

For a given partition $\lambda \vdash n$, any odd integer $x$, for any $t\geq 0$, call the number of times that $x\cdot 2^{t}$ appears in $\lambda \vdash n$  as $n^{\lambda}_{(x,t)}$.

Corresponding to this partition, for every odd integer $x$, construct the following matrix $M^{(x,\lambda)}$, where the $j'$th column from the left contains the description of $n^{\lambda}_{(x,j)}$: the $(i,j)$ th cell contains 0 if in the base 2 expansion of $n^{\lambda}_{(x,j)}$, the coefficient of $2^{i}$ is 0, and contains 1 otherwise. In other words, an entry $1$ in the $(i,j)$'th cell of $M^{(x,\lambda)}$ corresponds to $2^{i}$ many parts in the partition equal to $x\cdot 2^{j}$. For convenience, we label the rows and columns with indices $i,j$ each starting with 0. Henceforth we call the row indexed by $i$ as $\text{row}_{i}$, and the column indexed by $j$ as $\text{column}_{j}$.

\bigskip

For example, if for $\lambda \vdash n$, we have $n^{\lambda}_{(3,0)}=3$, $n^{\lambda}_{(3,1)}=6$ and $n^{\lambda}_{(3,2)}=3$, and $n^{\lambda}_{3,t}=0$ for all $t\geq 3$, then the matrix $M^{(3,\lambda)}$ looks like:

\begin{center}
\centering
\begin{tabularx}{0.4\textwidth}
 { | >{\raggedright\arraybackslash}X
  | >{\centering\arraybackslash}X
  | >{\centering\arraybackslash}X
  | >{\raggedleft\arraybackslash}X | }
 \hline
 $M^{(3,\lambda)}$ & 3 (j=0) & $3\cdot 2$ (j=1) & $3\cdot 2^{2}$ (j=2) \\
 \Xhline{4\arrayrulewidth}
 (i=0) & 1  & 0  & 1  \\
\hline
 (i=1)  & 1  & 1  & 1 \\
\hline
 (i=2) & 0  & 1  & 0 \\
\hline
\end{tabularx}
\end{center}

\bigskip


Consider for any integer $k\geq 0$ the ``diagonal" $D_{k}=\{(i,j): i+j=k \}$. When we permute the entries within any such fixed diagonal of any specific matrix corresponding to some odd number $x$, the contribution to the sum of the corresponding parts remains the same, but we get new a partition of $n$ in the process.

For the particular case above and the matrix $M^{(3,\lambda)}$, we could consider nontrivial permutations within the diagonal $D_1, D_2$ and $D_3$ to get different matrices, where we end up getting new partitions of $n$ in the process. However the single element diagonal $D_0$, and diagonals $D_t$ for $t\geq 4$ which are not shown in the figure for $t\geq 5$ which only contain $0$'s, remain invariant under any permutation.

For example, if we keep all diagonals except $D_2$ fixed, and consider the permutation within $D_2$ that interchanges the entries in cells $(2,0)$ and $(1,1)$, then we get a new partition $\lambda' \vdash n$ for which all the matrices $M^{(x,\lambda')}$ and $M^{(x,\lambda)}$ are identical when $x\neq 3$, and $M^{(3,\lambda')}$, written out below:
\bigskip

\begin{center}

\begin{tabularx}{0.4\textwidth}
 { | >{\raggedright\arraybackslash}X
  | >{\centering\arraybackslash}X
  | >{\centering\arraybackslash}X
  | >{\raggedleft\arraybackslash}X | }
 \hline
 $M^{(3,\lambda')}$ & 3 (j=0) & $3\cdot 2$ (j=1) & $3\cdot 2^{2}$ (j=2) \\
\hline
 (i=0) & 1  & 0  & 1  \\
\hline
 (i=1)  & 1  & 0  & 1 \\
\hline
 (i=2) & 1  & 1  & 0 \\
\hline
\end{tabularx}

\end{center}

For this partition, we have $n^{\lambda'}_{(3,0)}=2^0 +2^1+2^2=7$,  $n^{\lambda'}_{(3,1)}=2^2=4$, $n^{\lambda'}_{(3,2)}=2^{0}+2^{1}=3$.

\bigskip

First we briefly outline a proof of \cref{thm:oldone} that first appeared in \cite{Pathak}. For convenience of notation we deal with the $d=2$ case here, but the proof for general $d$ is identical. We show that the set of partitions $c'(n)$ of $n$ where at least one part is even, and the set of partitions $e'(n)$ where at least one part is repeated, have the same cardinality. Here we get a richer family of bijections that in the case of \cref{thm:conj3} whose proof follows later, where we fix the number of parts that are even and the number of parts that are repeated.

\begin{proof}[Outline of proof of \cref{thm:oldone}:] As stated before, without loss of generality we just consider the case $d=2$. Consider any partition in the set $c'(n)$. In each matrix of the form $M^{(x,\lambda)}$, within in each diagonal $D_k$ we permute all the elements so that the element in the cell $(k,0)$  goes to any other arbitrary cell within the same diagonal $D_k$, and the element in the cell $(0,k)$ is taken to the cell $(k,0)$. This gives us an element of $e'(n)$. The inverse map just takes the inverse permutations in each of the diagonals in each of the matrices of the form $M^{(x,\lambda)}$. It is easily verified that we have a bijection, whatever permutations we choose to take in each of these diagonals.
\end{proof}

We begin with the proof of \cref{thm:conj3}. By looking at these matrix constructions, the idea borrowed from \cite{Pathak}, the proof will be immediate.

\begin{proof}[Proof of \cref{thm:conj3}:]
Consider the set $A_{k}(n)$ of partitions that have exactly $k$ distinct even parts. The binary description of the number of appearances of each of these distinct even parts are distributed among $k$ different columns in total, distributed among one or more of the matrices constructed above. Every other column with index $j\geq 1$ in any of the matrices has all elements 0. Thus we are only concerned with columns with indices $j\geq 1$ whereas $\text{column}_0$ (corresponding to $j=0$) may have any number of entries that are $1$.

 For each of the filled columns, push each entry of the specific column diagonally one place down and left $(i.e \ (i,j)\to (i+1,j-1))$, and for each $i\geq 1$, take the $(i,0)$ entry (entries in the $\text{column}_0$) to the $(0,i)$ entry (entries in $\text{row}_{0}$), and keep the $(0,0)$ entry constant. It is clear that this is a bijection, taking an element of $A_{k}(n)$ to an element of $B_{k}(n)$, which is the set of partitions where exactly $k$ different parts are repeated. The inverse map takes, for each $i\geq 1$, the entries of cell $(i,j)$ to the $(i-1, j+1)$'th cell, keeps the $(0,0)$'th entry constant, and for each $j\geq 1$, takes the element in $(0,j)$'th cell to the $(j,0)$ cell.
\end{proof}

To illustrate, we consider a generic partition  $\lambda\vdash n$ belonging to the set $A_{j}(n)$ which has several matrices of the form $M^{(x_{i},\lambda)}$ where $i$ ranges in some finite set.

 On the left, we have a generic matrix corresponding to the partition $\lambda\vdash n$ belonging to $A_{j}(n)$ and on the right, we have the matrix of the corresponding partition $\lambda' \vdash n$ belonging to $B_{j}(n)$ that is obtained through the bijection. The entries labelled as $a_{ij}$ can take the values $0$ or $1$, where $i,j\geq 0$ and the entries are eventually zero for large enough $i.j$ all 0.

\begin{table}
\centering
\begin{tabularx}{0.4\textwidth}
 { | >{\raggedright\arraybackslash}X
  | >{\centering\arraybackslash}X
  | >{\centering\arraybackslash}X
  | >{\raggedleft\arraybackslash}X | }
 \hline
 $M^{(x,\lambda)}$ & $x$ (j=0) & $x\cdot 2$ (j=1) & $x\cdot 2^{2}$ (j=2) \\
 \Xhline{4\arrayrulewidth}
 (i=0) & $a_{00}$  & $a_{01}$  & $a_{02}$  \\
\hline
 (i=1)   & $a_{10}$  & $a_{11}$  & $a_{12}$ \\
\hline
 (i=2)  & $a_{20}$  & $a_{21}$ & $a_{22}$ \\
\cdashline{1-4}
\end{tabularx}
$\longleftrightarrow$
\begin{tabularx}{0.4\textwidth}
 { | >{\raggedright\arraybackslash}X
  | >{\centering\arraybackslash}X
  | >{\centering\arraybackslash}X
  | >{\raggedleft\arraybackslash}X | }
 \hline
 $M^{(x,\lambda')}$ & $x$ (j=0) & $x\cdot 2$ (j=1) & $x\cdot 2^{2}$ (j=2) \\
 \Xhline{4\arrayrulewidth}
 (i=0) & $a_{00}$  & $a_{10}$  & $a_{20}$  \\
\hline
 (i=1)  & $a_{01}$  & $a_{02}$  & $a_{03}$ \\
\hline
 (i=2) & $a_{11}$  & $a_{12}$  & $a_{13}$ \\
 \cdashline{1-4}
\end{tabularx}
\caption{On the left, the matrix $M^{(x,\lambda)}$ of the partition belonging to $A_{j}(n)$, and the corresponding matrix of the corresponding partition belonging to $B_{j}(n)$, obtained through the bijection of \cref{thm:conj3}. Only the matrix entries shown are possibly non-zero and the other entries for higher values of $i,j$ are all zero.}

\end{table}

Building on the above argument, we give a broader family of bijections where the partitions satisfy the more restrictive condition of \cref{thm:new1}.

\begin{proof}[Proof of \cref{thm:new1}:]
The argument here is a generalization of that of the previous proof. Consider any partition in $A_{k}(n)$, and any arbitrary matrix corresponding to an arbitrary odd number $x$.

We transfer the element in the cell $(i,j)$ to the cell $(i+p,j-p)$, when $j\geq p$. In the $p \times p$ block in the top left hand corner, we can permute the elements so that each element in any cell $(a,b)$ only moves along the diagonal $D_{(a+b)}$ containing it, similar to the case in the proof in \cite{Pathak}. Furthermore, in this top left block, the permutations can be different in different diagonals, and we would still get a valid bijection.\footnote{In this case except the diagonal $D_{p-1}$, all other diagonals are not completely within this $p\times p$ block and we are only restricting to the portion of the diagonal within this block.} It remains to permute the elements in the set of cells $\{(i,j): 0\leq j\leq p-1, i> p-1\}$ to the set of cells $\{(i,j): j> p-1, 0\leq i \leq p-1\}$. This can be done in several different ways, and each way gives a valid bijection. For each integer $m\geq 2$ to take the blocks $\{(i,j):mp\leq i \leq (m+1)p-1, 0\leq j\leq p-1 \}$, permute the elements within each diagonal within this block, and transpose \big(i.e. taking the entries in the cells $(i,j) \leftrightarrow (j,i)$ \big) it to the set $\{(i,j): mp\leq j \leq (m+1)p-1, 0 \leq i\leq p-1 \} $. Just as in the previous proof, the inverse maps are also obvious for each of this family of bijections. For each $i\geq p$, the inverse bijection takes the element in cell $(i,j)$ to the cell $(i-p,j+p)$, while in the top left $p\times p$ block, along each diagonal we take the inverse permutation. The elements of the set of cells $\{(i,j): 0\leq i \leq p-1, j> p-1 \}$ are transposed (i.e. taking the element of the cell $(i,j)\leftrightarrow (j,i)$) to the set of cells  $\{(i,j): 0\leq j <p-1, i\geq p-1 \}$, and then in each of the blocks $\{(i,j):mp\leq j \leq (m+1)p-1, 0\leq i \leq p-1 \}$ we apply the inverse permutation. \footnote{To construct the inverse map, in effect we could also carry out a reverse permutation first and then transpose.}
\end{proof}

An illustration of an instance of this bijection follows closely the illustration of the previous bijection, except for the freedom to permute within the diagonals of the blocks where either of the $i$ or $j$ coordinates are restricted within the range $[0,p-1]$, as explained in the proof.

\begin{table}
\centering
\begin{tabular}{|c|c|c|c|}
\hline
$M^{(x,\lambda)}$ & $(j=0,\dots, p-1)$ & $(j=p,\dots,2p-1)$ & $(j=2p,\dots,3p-1)$ \\
\Xhline{4\arrayrulewidth}
 $(i=0,\dots, p-1)$ & $M_{00}$  & $M_{01}$  & $M_{02}$  \\
\hline
 $(i=p,\dots, 2p-1)$ & $M_{10}$  & $M_{11}$  & $M_{12}$ \\
\hline
 $(i=2p, \dots, 3p-1)$ & $M_{20}$  & $M_{21}$  & $M_{23}$ \\
 \cdashline{1-4}
\end{tabular}

\bigskip

$\updownarrow$

\bigskip

\begin{tabular}{|c|c|c|c|}
\hline
 $M^{(x,\lambda')}$ & $(j=0,\dots, p-1)$ & $(j=p,\dots, 2p-1)$ & $j=(2p,\dots,3p-1)$ \\
 \Xhline{4\arrayrulewidth}
 $(i=0,\dots, p-1)$ & $\overline{M}_{00}$  & $\overline{M}_{10}$  & $\overline{M}_{20}$  \\
\hline
 $(i=p,\dots, 2p-1)$  & $M_{01}$  & $M_{02}$  & $M_{03}$ \\
\hline
 $(i=2p, \dots, 3p-1)$ & $M_{11}$  & $M_{12}$  & $M_{13}$ \\
 \cdashline{1-4}
\end{tabular}

\bigskip

\caption{Above, the matrix $M^{(x,\lambda)}$ of the partition $\lambda\vdash n$ belonging to $A_{k}(n)$, and the corresponding matrix of the corresponding partition $\lambda'\vdash n$ belonging to $B_{k}(n)$, obtained through the bijection of \cref{thm:new1}. Each of the entries above are $p\times p$ matrices, where $\overline{M}_{00}$ is obtained from $M_{00}$ by arbitrary permutations on the diagonals. The matrices $\overline{M}_{10}$ and $\overline{M}_{20}$ are respectively found from $M_{10}$, $M_{20}$ through arbitrary permutations along its diagonals and then transposing.}

\bigskip

\end{table}

Note that any of the general family of bijections for \cref{thm:new1} where $p\geq 2$ stated above are not useful for proving \cref{thm:conj3} where $p=1$, since in the latter case any $1\times 1$ block consists of a single element and thus permutations within such a block is trivial. Also, instead of considering powers of $2$, we could easily deal with the powers of any arbitrary integer $d$ and the statements of the above theorems are suitably modified. This analogous generalization for \cref{thm:conj3} was carried out in Fu and Tang's paper \cite{Fu}, and the similar generalization of \cref{thm:new1} is also immediate.

\section{Combinatorial proof of Theorem 2 and generalizations.}

We first give a direct short proof of \cref{thm:conj2} using broadly the same methods as before. After that we state and prove two different generalizations of this result.

\begin{proof} [Proof of  \cref{thm:conj2}:]
It is clear that the number that is the difference between the total number of parts in the distinct partitions of $n$ and the total number of distinct parts in the odd partitions of $n$, can be split up as a sum of smaller parts as we outline now.

Call the set of partitions of $n$ with all distinct parts $D_{n}$, and the set of partitions of $n$ with all odd parts as $O_{n}$. Given any partition $\lambda\vdash n$ belonging in $D_{n}$, for some fixed odd number $x$, the matrix $M^{(x,\lambda)}$ has some number of $1$'s in the first row, and the matrix $M^{(x,\lambda)}$ otherwise contains all zeros. Call this number of $1$'s in the first row as $d_{x}$. In $O_{n}$, the entries of the previous matrix are essentially flipped to come to the first column of $M^{(x,\mu)}$ where $\mu$ is the partition obtained from $\lambda$ through Glaisher's bijection. Each such matrix $M^{(x,\mu)}$ gives a count of 1 to  the number of distinct parts in any partition in $O_{n}$. Thus, per fixed matrix of the form $M^{(x,\lambda)}$, to the difference of the number of parts in the distinct partitions of $n$ and the number of distinct parts in the odd partitions of $n$, we get a contribution $(d_{x}-1)$. When we sum such contributions over all the possible matrices of the form $M^{x,\lambda}$ as $x, \lambda$ are varied, we get the total difference.


Hence, given partition of $\lambda\vdash n$ into distinct parts, and odd integer $x$, its corresponding matrix $M^{(x,\lambda)}$, consider only the cells $(0,j_1),(0,j_2),\dots,(0,j_m)$ to have $1$'s \footnote{$m$ equals the $d_x$ introduced in the previous paragraph. }, for some positive integer $m$, and between $(0,j_{a})$ and $(0,j_{a+1})$ let there be $t_a\geq 0$ empty cells, where $1\leq a\leq m-1$ are integers.

We show that corresponding to the $M^{(x,\lambda)}$ above, we can get $(m-1)$ many partitions in each of which there is exactly one part that appears three times and all others appear once, and when we add the corresponding contributions from $M^{(x,\lambda)}$ as $\lambda, x$ are varied,  we get the total number of partitions in $a_{1}(n)$.

Indeed, in the above matrix $M^{(x,\lambda)}$, consider some $l$ with $1\leq l \leq m-1$, and a new partition of $n$ where in the first row all the filled and unfilled cells prior to the $(0,j_{l})$ cell remain unchanged, and the cells $(0,j_{l})$ and $(1,j_{l})$ are filled with 1, and all the originally empty cells of the first row in between  $(0,j_{l})$ and $(0,j_{l+1})$ are filled with 1's and the originally filled cell $(0,j_{l+1})$ now contains 0. It is simple to check that the sum contributed by this configuration is the same as the sum contributed by the original configuration.\footnote{this is the only way here to keep the sum constant and to have a case where there is one part appearing three times and all other parts appearing once}.

Thus we have shown that corresponding to this matrix configuration, we can get $(m-1)$ many partitions in which there's exactly one part that appears three times and all others appear once. Also it is clear that summing on all the  distinct matrices, we get the exact count for the number of partitions in this set of partitions $a_{1}(n)$. Thus the proof is complete.
\end{proof}

Now we generalize the statement of \cref{thm:conj2} to give the proof of \cref{thm:new2}.

\begin{proof} [Proof of \cref{thm:new2}:]
The proof of this theorem essentially uses the same argument as the previous proof. Again for a partition $\lambda\vdash n$ belonging in the set $G(n)$ and a specific odd $x$, if there are $m$ many appearances of $1$ in the first row in the matrix $M^{(x,\lambda)}$, say in the cells $(0,j_1), (0,j_2),(0,j_3),..,(0,j_m)$ where now there is at least one empty cell between each of these filled cells, then in order to get distinct elements of the set of partitions with only one part repeated five times, all we do is choose a specific $1\leq a \leq m-1$, put $1$'s in the cells $(0,j_a)$ and $(2,j_a)$, keep the cell $(0,j_{a}+1)$ empty and put $1$'s in all the cells from $(0,j_{a}+2)$ to $(0,j_{a+1}-1)$ and make $(0,j_{a+1})$ empty. In that process we get the new partition, and the theorem follows.
\end{proof}

In Yang's paper \cite{Yang}, there is a generalization of \cref{thm:conj2} in the Glaisher like setting with an arbitrary integer $d$, in Theorem 1.7 of that paper, where the set of partitions is considered where exactly one part appears more than $d$ times and less than $2d$ times whereas all other parts appear less than $d$ times. The proof of this would essentially follow our proof technique of \cref{thm:conj2} above. Combining this with our \cref{thm:new2} above, we can formulate other interesting statements in the Glaisher setting.

\section{Combinatorial proof of Theorem 1.}

In this last section, we give quick combinatorial arguments for \cref{thm:conj1}.

\begin{proof}[First proof of \cref{thm:conj1}:]
In this case, consider the number of times $n_x$, that a particular odd number $x$ appears in one of the odd partitions of $n$, and write it in the  base 2 expansion: $n_x=\gamma_0 +2\cdot \gamma_1 +2^{2}\cdot \gamma_2+\dots+2^{m}\cdot \gamma_m$, with $\gamma_m=1$, where each of the other $\gamma_i$ are either 0 or 1. By an argument similar to that used in the previous proof, it should be clear that the difference between the total number of parts in the odd partitions of $n$ and the number of parts in the distinct partitions of $n$, can be broken up into summands per matrix $M^{(x,\lambda)}$ as $\lambda\vdash n$ is varied within the odd partitions, $x$ is varied over the odd integers. The summands are of the form $(\gamma_0 +2\cdot \gamma_1 +2^{2}\cdot \gamma_2+\dots+2^{m}\cdot \gamma_m)-(\gamma_0 + \gamma_1 + \gamma_2+\dots+ \gamma_m)= (2-1)\cdot \gamma_1 +(2^{2}-1)\cdot \gamma_2+\dots+(2^{m}-1)\cdot \gamma_m$.

Now we show that per matrix $M^{(x,\lambda)}$ above, we will find $(2-1)\cdot \gamma_1 +(2^{2}-1)\cdot \gamma_2+\dots+(2^{m}-1)\cdot \gamma_m$ distinct partitions that contain exactly one even part, where that even part may be repeated more than once.
Consider the odd partitions of $n$ where $x$ appears $n_x$ number of times. The matrix has $1$'s in only the first column, and all 0's in all other columns. We wish to introduce exactly one even part here; i.e. introduce some number of $1$'s into exactly one more column, while keeping the sum constant.

For the above fixed matrix, consider the set of indices  $\{i: \gamma_i=1\}$, and label these coefficients that are $1$, as $\gamma_{i_1},\gamma_{i_2},\gamma_{i_3},\dots \gamma_{i_{k}}$. Thus in the first column, only the cells $(i_1,0),(i_2,0),\dots$ $,(i_k,0)$ contain 1 and the rest contain 0. Consider first $i_1$.  In the second column, one can put any number of $1$'s in the cells $(0,1),(1,1),\dots,(i_1 -2,1)$; then we will put the binary description of the difference of this new obtained number from $2^{i} $, in the first column, keeping the entry in the cell $(i_1 -1,1)$ equal to 0 and the coefficients of $2^{j}$ in the first column for $j\geq i_1 +1$ same as before. Otherwise, we could also put a $1$ in the $(i_1 -1,1)$ cell, with all entries in the cells $(0,1),\dots,(i_1-2 ,1)$ being $0$, while putting $0$ in the $(i,0)$'th cell and keeping everything else in the first column same as before.  In total, we thus have $(2^{i_1-1}-1)+1=2^{i_1 -1}$ choices to fill up the second column, and keeping the coefficients of  $2^{j}$ for $ j\geq i+1$ unchanged in the first column.

Now consider $i_2$ and the second column. In this case, for all $j\geq i_{2}+1$ we keep the same entry as before in the $(j,0)$ cells. We count the cases where in the second column we can put any number of $1$'s in the first $i_1 -1$ cells from the top, while keeping the coefficient of $2^{i_1}$ equal to 1; ( i.e. the cell $(i_1 -1,1)$ having the entry 1), the coefficient of $2^{i_2}$ equal to 0; (i.e. the cell $(i_2-1,1)$ having the entry 0). Next we add to this previous count, the count of cases where we have at least one entry being 1 among the cells of the second column with row numbers greater than $i_1 -1$ and less than $i_{2}-1$. Its obvious that we have thus $2^{i_2 -1} -1$ many new cases here. Finally consider the case where we have everything below the $(i_2 -1,1)$ cell in the second column having an entry 0 and just the $(i_2 -1,1)$ cell having the entry $1$. Thus in total we have $2^{i_2-1}$ cases from here.

Continuing this algorithm for $i_k, k\geq 3$, it should be clear that for the second column, corresponding to every $i_{k}$ we would have $2^{i_k -1}$ cases to count.

Now if we wanted to fill the $m$'th column with $m\geq 3$ instead, we would have for any given $k$,  $2^{i_k-(m-1)}$ choices corresponding to $i_{k}$. On the other hand, for each fixed $i_{k}$, summing over all the columns, we get $2^{i_k -1} +2^{i_{k}-2}+\dots +1=2^{i_{k}}-1$ choices.

Thus finally we find (trivially including the $\gamma_i$ that are also $0$) the number $(2-1)\cdot \gamma_1 +(2^{2}-1)\cdot \gamma_2+\dots +(2^{m}-1)\cdot \gamma_m$.
\end{proof}

Below we give a quick outline of a slight variation of the argument of the above proof to arrive at the same count.

\begin{proof}[Second proof of \cref{thm:conj1}:]

Consider again the number $n_\alpha=\gamma_0 +2\cdot \gamma_1 +2^2 \cdot \gamma_{2}+\dots+2^{m}\cdot \gamma_{m}$, from the previous proof, which is the number of times the odd number $\alpha$ appears in one of the odd partitions of $n$. Assume that we want to break this up into two parts, with $\alpha$ appearing some number of times, and $2^{j}\cdot \alpha$ (for any $j\geq 1$) appearing some other number of times, and these together summing up to give $\alpha\cdot n_\alpha$. Rewrite $n_{\alpha}$ as $n_\alpha= (\gamma_0 +2\gamma_1 +\dots+\gamma_{j-1}2^{j-1})+2^{j}(\gamma_j +2\cdot\gamma_{j+1}+\dots )$. It is clear that we have $(\gamma_j +2\gamma_{j+1}+\dots +2^{m-j}\gamma_{m})$ many choices, which is the number in the second bracket above. Thus by an easy double counting, for any fixed coefficient $\gamma_{t}$, as we vary the $j$'s above, we have $1+2+2^2 +\dots+2^{t-1}=2^{t}-1$  appearances of $\gamma_t$. Thus as before we have the total number $(2-1)\cdot \gamma_1 +(2^2 -1)\cdot \gamma_2 +\dots+(2^{m}-1)\cdot \gamma_m$.
\end{proof}


\bigskip



\end{document}